\documentclass[11pt]{article}

\usepackage{amsmath,amssymb,amsthm}
\usepackage{mathrsfs}
\usepackage{mathtools}
\usepackage{geometry}
\usepackage{breqn}
\usepackage[numbers,sort&compress]{natbib}
\usepackage{stmaryrd}
\usepackage[hidelinks]{hyperref}
\usepackage{cleveref}
\usepackage{authblk}
\usepackage{comment}
\DeclarePairedDelimiter\floor{\lfloor}{\rfloor}
\DeclarePairedDelimiter{\ground}{\llbracket}{\rrbracket}
\DeclarePairedDelimiter\croc{\langle}{\rangle}

\theoremstyle{definition}

\theoremstyle{remark}

\newcommand{\ii}{\mathrm i}
\newcommand{\dd}{\mathrm d}
\newcommand{\ee}{\operatorname{e}}
\newcommand{\ev}{\operatorname{ev}_{0}}
\newcommand{\RE}{\operatorname{Re}}

\makeatletter
\newcommand*{\defeq}{\mathrel{\rlap{%
  \raisebox{0.3ex}{$\m@th\cdot$}}%
  \raisebox{-0.3ex}{$\m@th\cdot$}}%
  =}
\makeatother

\title{Umbral Hermite Representations of Canonical Catastrophe Integrals\\
\vskip 0.25cm
\large Formal Operators, Analytic Pairings and Scalar Partition Kernels}
\author[1]{Giuseppe Dattoli\footnote{E-mail: pinodattoli@libero.it (Giuseppe Dattoli)}}
\author[1]{Roberto Ricci\footnote{E-mail: roberto.ricci@enea.it (Roberto Ricci)}}
\affil[1]{ENEA, Nuclear Department, Frascati Research Center, Via E. Fermi 45, 00044 Frascati (Rome), Italy}

\date{}

\begin{document}

\maketitle

\begin{abstract}
Higher-order multivariable Hermite--Kamp\'e de F\'eriet polynomials
provide a unified construction of umbral representations for canonical
catastrophe integrals. Expanding the unfolding exponential in these
polynomials and integrating against the leading catastrophe kernel
gives convergent moment expansions on suitable decay contours.
A universal formal umbral operator then expresses the resulting
fields through explicit moment ground states. We develop this
construction for the fold, cusp and swallowtail, and relate the
Euclidean Pearcey field to quartic scalar partition functions.

The distinction between formal and analytic umbral frameworks becomes
essential when integration is performed before umbral evaluation.
For Euclidean Pearcey, the integrated linearized kernel is a rational
function whose action is not specified by the original formal rules.
The Analytic Umbral Framework supplies a convergent Mellin pairing
for this kernel and recovers the independently obtained bivariate
Pearcey representation. This establishes, in this example, an
exchange of integration and analytic evaluation and identifies a
specific extension beyond the formal framework.
\end{abstract}

\section{Introduction}
\label{sec:introduction}

Canonical catastrophe integrals describe universal diffraction
patterns associated with coalescing stationary points of an
oscillatory phase
\cite{Thom1975,Arnold1992,BerryUpstill1980}.
The fold, cusp and swallowtail correspond to the $A_2$, $A_3$
and $A_4$ canonical integrals, respectively; the first two are
expressed through the Airy and Pearcey functions
\cite{Pearcey1946,ValleeSoares2004}.
Related polynomial integrals also occur in statistical mechanics
and field theory. In particular, the Euclidean quartic integral
provides a common analytic description of the Pearcey field and
zero-dimensional scalar partition functions.

The objective of this paper is to construct umbral representations
of canonical catastrophe fields, starting from their standard
oscillatory definitions. The organizing objects are the higher-order
multivariable Hermite--Kamp\'e de F\'eriet (HKdF) polynomials.
Their generating function expands the lower-degree, parameter-dependent
part of the phase, while the leading exponential is retained as
the integration kernel. Integration then replaces each power of
the integration variable by a canonical contour moment.
The result is a Hermite expansion of the catastrophe field.
On suitable decay contours, the exchange of summation and
integration is justified locally uniformly in the unfolding
parameters.

The umbral formulation packages this construction into compact
operator expressions. In its traditional form, higher-order
Hermite umbral numbers are represented by operators $\hat h_k$
acting on symbolic vacua $\varphi_0^{(k)}$
\cite{LicciardiDattoli2022,DattoliKhanAhmad2025}.
We recast these evaluation rules using a universal formal umbral
operator $\mathfrak u$ and explicit order-dependent ground state functions
$\varphi_k(t)$ \cite{Ricci2026}.
After canonical integration, the moment sequence supplies a new
ground state for the field itself. One may retain a single umbra,
with polynomial dependence in its exponent, or introduce one
commuting umbra for each unfolding parameter and encode their
coupling in a multivariate ground state.

This formal construction does not replace the geometric content
of catastrophe theory. The phase normal form, unfolding parameters
and integration cycle remain part of the definition of the
canonical field. Contour choice is also essential analytically:
a decaying contour makes the moment construction well defined,
whereas high polynomial moments on the undeformed real axis
need not exist as ordinary improper integrals. We keep this
distinction explicit in the fold, cusp and swallowtail examples.

The Analytic Umbral Framework (AUF) realizes evaluation through
Mellin--Barnes pairings \cite{RicciAUF2026}. It introduces analytic
requirements that are absent from the formal rules: the
transmutation kernel must admit an appropriate Mellin transform,
and the resulting pairing must have compatible branches and
sufficient convergence. Convergence of the Hermite series alone
does not guarantee these properties for an arbitrary Mellin
representation. The Euclidean Pearcey case allows them to be
checked explicitly.

Its particular significance is revealed by reversing the order
of integration and evaluation. In the natural order, the formal
umbral framework (FUF) reconstructs the polynomial exponential
before ordinary contour integration, and the calculation is exact.
Integrating the linearized expression first instead suggests a
rational function of a sum of umbral quantities. The original
formal rules neither define the required operator integral nor
determine the evaluation of this rational expression.

For Euclidean Pearcey, the AUF supplies the missing analytic
extension. We integrate the scalar transmutation kernel on
specified paths, calculate its Mellin transform and pair the
result with the monomial ground states. The resulting convergent
integral reproduces the independently established bivariate
Pearcey representation. Thus the two orders agree at the level
of analytic pairings. This provides a precise realization of
the operational calculation without assuming that an
unevaluated formal operator integral already exists.

The paper is organized as follows.
Section~\ref{sec:HKdF} recalls the higher-order multivariable HKdF polynomials. Section~\ref{sec:umbral} introduces
evaluation rules for these polynomials, their traditional and universal formal realizations and the
AUF pairing.
Section~\ref{sec:canonical} specifies the canonical $A_k$
integrals, their contours and their Hermite moment expansions.
Section~\ref{sec:representations} develops the resulting umbral
representations, with the fold, cusp and swallowtail as explicit
examples and quartic partition functions as an application.
Section~\ref{sec:reversed} establishes the reversed-order
AUF calculation for Euclidean Pearcey.
Section~\ref{sec:conclusions} summarizes the results and prospects.

\section{Higher-order multivariable Hermite polynomials}
\label{sec:HKdF}

The two-variable Hermite-Kamp{\'e} de F{\'e}riet (H-KdF) polynomials of degree $n$ and order $N \geq 2$ are defined by
\begin{dmath}\label{eq:Hermite-Kampe}
	H_n^{[N]}(x, y) \defeq n!\sum_{r=0}^{\floor{\frac{n}{N}}}
	\frac{
	x^{n-Nr} y^r
	}{
	(n-Nr)! r!
	}
\end{dmath}.
For $N > 2$, they represent the higher-order extension of the standard two-variable Hermite polynomials $H_n(x,y) \equiv H_n^{[2]}(x, y)$.
Their generating function is
\begin{dmath}\label{eq:HK-generating-function}
	\sum_{n=0}^\infty 
	H_n^{[N]}(x, y)
	\frac{t^n}{n!}
	=
	\ee^{xt+yt^N}
\end{dmath}.
They also possess the operational representation \cite{LicciardiDattoli2022,DattoliKhanAhmad2025}
\begin{dmath}\label{eq:HK-operational}
	H_n^{[N]}(x, y) 
	=
	\ee^{y\partial_x^N} x^n
\end{dmath}.
These formulas admit an immediate multivariable extension \cite{DattoliLicciardiSabia2025,DattoliKhanAhmad2025}. In particular, when the number of variables is $N$, 
\begin{dmath}\label{eq:Hermite-Kampe-multivariable}
	H_n^{[N]}(\boldsymbol x) \equiv 
	H_{n}^{[N]}(x_{1},x_{2},\dots ,x_{N}) = 
	n!\sum_{\substack{
	r_{2},\dots ,r_{N} \geq 0\\
	2r_2+\cdots+Nr_N\leq n}}
	\frac{x_{1}^{n-2r_{2}-\dots -N r_{N}}\cdot x_{2}^{r_{2}}\cdots x_{N}^{r_{N}}}{(n-2r_{2}-\dots -N r_{N})!\cdot r_{2}!\cdot r_{3}!\cdots r_{N}!}
	=
	n! \sum_{\substack{r_2, \dots, r_N \geq 0 \\ 2r_2 + \dots + Nr_N \leq n}} \frac{x_1^{r_1} \cdot x_2^{r_2} \cdots x_N^{r_N}}{r_1! \cdot r_2! \cdots r_N!}
	=
	n!
	\sum_{
	\substack{r_1,\ldots,r_N\geq0\\
	r_1+2r_2+\cdots+Nr_N=n}}
	\prod_{k=1}^{N}
	\frac{x_k^{r_k}}{r_k!}
\end{dmath}.
The corresponding generating function is
\begin{dmath}\label{eq:multivariable-generating-function}
\sum_{n=0}^\infty
H_n^{[N]}(\boldsymbol x)
\frac{t^n}{n!}
=
\exp\!\left(
\sum_{k=1}^N x_k t^k
\right)
\end{dmath}.
Raising and lowering rules follow directly from \cref{eq:multivariable-generating-function}:
\begin{dmath}
	\frac{\partial}{\partial x_1} H_n^{[N]}(\boldsymbol x) = n H_n^{[N]}(\boldsymbol x)
\end{dmath},
and, more generally,

\begin{dmath}\label{eq:multivariable-lowering}
	\frac{\partial}{\partial x_j}
	H_n^{[N]}(\boldsymbol x)
	= {
	\frac{n!}{(n-j)!}
	H_{n-j}^{[N]}(\boldsymbol x),
	\qquad
	1 \leq j \leq n
	}
\end{dmath}.
The operational representation \cref{eq:HK-operational} generalizes to
\begin{dmath}
	H_n^{[N]}(\boldsymbol x) =
	\exp
	\left(
	\sum_{k=2}^N
	x_k\partial_{x_1}^k
	\right)
	x_1^n
\end{dmath}.

The $N$-variable H-KdF polynomials satisfy the recurrence relation
\begin{dmath}\label{eq:multivariable-recurrence}
	H_{n+1}^{[N]}(\boldsymbol x)
	=
	x_1H_n^{[N]}(\boldsymbol x)
	+
	\sum_{j=2}^{N}
	jx_j
	\frac{n!}{(n-j+1)!}
	H_{n-j+1}^{[N]}(\boldsymbol x)
\end{dmath}.

The higher-order multivariable HKdF polynomials are rescaled
specializations of the complete exponential Bell polynomials \cite{Comtet1974}:
\begin{dmath*}
	H_n^{[N]}(\boldsymbol x)
	= {
	B_n(y_1,\ldots,y_n),
	\qquad
	y_j=
	\begin{cases}
		j!\,x_j, & j\leq N,\\
		0, & j>N
	\end{cases}}
\end{dmath*}.
Thus the two families share the same generating-function structure. The HKdF notation emphasizes the finite-degree polynomial exponent and its operational interpretation.
This identification locates the HKdF family within classical combinatorics. The specific contribution of the present construction is the multicomponent umbral realization adapted to canonical catastrophe phases, rather than the polynomials considered in isolation.

\section{Formal and analytic umbral representations}\label{sec:umbral}

The H-KdF polynomials admit a formal umbral representation that enables to greatly simplify the previous formulas.

For \(2\leq k\leq N\), introduce the formal operators \(\{\hat h_k\}\) and the corresponding symbolic vacua \(\{\varphi_0^{(k)}\}\) through
\begin{dmath}\label{eq:formal-moments}
\hat h_k^n\varphi_0^{(k)}
=
\begin{cases}
\dfrac{n!}{(n/k)!}, & {n=kq,\quad q\in\mathbb N_0},\\
0, & \text{otherwise},
\end{cases}
\end{dmath}.

Using \cref{eq:formal-moments}, \cref{eq:Hermite-Kampe} can be written as
\begin{dmath}\label{eq:HK-umbral}
	H_n^{[N]}(x, y) 
	=
	\left(
	x + y^{1/N} \hat h_N
	\right)^n
	\varphi_0^{(N)}
\end{dmath}.

The umbral construction can be generalized to the multidimensional case by assuming that the $N-1$ operators are independent and commuting. Introduce the factorized vacuum
\begin{dmath*}
	\boldsymbol\varphi_0
	=
	\prod_{k=2}^{N}\varphi_0^{(k)}
\end{dmath*}
and the distinct linear combination
\begin{dmath}\label{eq:formal-polynomial-component}
	\widehat{H}_N(\boldsymbol x)
	\defeq
	x_1
	+
	\sum_{k=2}^{N}(x_k)^{1/k}\hat h_k
\end{dmath}.
Then
\begin{dmath}\label{eq:HK-multidimensional-umbral}
	H_n^{[N]}(\boldsymbol x) = {
	\left(
	x_1 + x_2^{1/2} \hat h_2 + \cdots + x_N^{1/N} \hat h_N
	\right)^n
	\boldsymbol\varphi_0 =
	\left(\widehat{H}_N(\boldsymbol x)\right)^n \boldsymbol\varphi_0
	}
\end{dmath}.
This expression is the defining umbral realization used in the following analysis. 
The polynomials $H_n^{[N]}(\boldsymbol x)$ enable to separate the algebra of the unfolding parameters from the moments of the leading catastrophe kernel, leading to explicit umbral representations of the canonical catastrophe integrals.

Consider the formal operator $\ee^{s\hat h_k}$. Expanding the exponential gives
\begin{dmath}\label{eq:fundamental-formal}
	\ee^{s\hat h_k}
	\varphi_0^{(k)}
	=
	{
	\sum_{q=0}^{\infty}
	\frac{s^{kq}}{q!}
	=
	\ee^{s^k}
	}
\end{dmath}.
For a polynomial
\begin{dmath*}
	P_N(s; \boldsymbol c)
	=
	c_1s
	+
	\sum_{k=2}^{N}c_ks^k
\end{dmath*},
we obtain
\begin{dmath}\label{eq:formal-linearization}
	\ee^{s\widehat{H}_N(\boldsymbol c)}
	\boldsymbol\varphi_0
	=
	\ee^{P_N(s; \boldsymbol c)}
\end{dmath}.
Although the individual roots require branch choices, the evaluated expression is
branch independent: changing \((c_k)^{1/k}\) by a \(k\)-th root of unity does
not affect the only surviving powers \(c_k^{q}\).

\subsection{Universal umbral operator and explicit ground states}

Equation~\eqref{eq:formal-moments} is operationally efficient but leaves
both \(\hat h_k\) and \(\varphi_0^{(k)}\) unspecified. The same rule can be
written in terms of a universal umbral operator $\mathfrak u$ and an explicit family
of ground-state functions.

Let \(\ev f=f(0)\). For a function \(f\) and a complex number \(\alpha\)
in its domain, define
\begin{dmath}\label{eq:universal-operator}
\mathfrak u^\alpha
\ground{f(t)}
\defeq
{
\left(
\ev\circ\exp(\alpha\partial_t)
\right)
f(t)
=
f(\alpha)
}
\end{dmath}.
The operator \(\mathfrak u\) is universal, and its definition can be made rigorous in the context of the theory of formal power series \cite{Ricci2026}.

For the problem considered in this work, define the congruence selector
\begin{dmath}\label{eq:selector}
\Pi_k(t_k)
=
\exp\!\left[
-\frac{\pi\ii(k-1)}{k}t_k
\right]
\frac{\sin(\pi t_k)}{k\sin(\pi t_k/k)}
\end{dmath},
with removable singularities filled by continuity, and set
\begin{dmath}\label{eq:explicit-ground-state}
\varphi_k(t_k)
= {
\frac{\Gamma(1+t_k)}{\Gamma(1+t_k/k)}
\Pi_k(t_k),
\qquad
2 \leq k \leq N
}
\end{dmath}.
Each ground state function satisfies the finite-difference equation
\begin{dmath}\label{eq:k-step-recurrence}
\varphi_k(t_k+k)
= 
k
\prod_{j=1}^{k-1}(t_k+j)
\varphi_k(t_k)
\end{dmath}.
The recurrence propagates all congruence classes; the selector
\(\Pi_k\) fixes the one required by the Taylor support.

At nonnegative integers,
\begin{dmath}\label{eq:ground-state-lattice}
\varphi_k(n)
=
\begin{cases}
\dfrac{n!}{(n/k)!}, & {n=kq,\quad q\in\mathbb N_0}, \\
0, & \text{otherwise}.
\end{cases}
\end{dmath}
Hence, by denoting with $\mathfrak u_k$ the universal operator acting on functions of $t_k$, we get
\begin{dmath}\label{eq:universal-equivalence}
\mathfrak u_k^n
\ground{\varphi_k(t_k)}
=
{
\varphi_k(n)
=
\hat h_k^n\varphi_0^{(k)}
}
\end{dmath},
and
\begin{dmath}\label{eq:universal-exponential}
\ee^{
s\mathfrak u_k
}
\ground{\varphi_k(t_k)}
=
\ee^{s^k}
\end{dmath}.

For a factorized polynomial exponential, we assume that the $N-1$ operators
\(\{\mathfrak u_k\}_{k=2}^N\) are independent and commuting. Introducing the factorized ground state function
\begin{dmath*}
	\varphi(\boldsymbol  t) \defeq \prod_{k=2}^{N}\varphi_k(t_k)
\end{dmath*}
and the distinct linear combination
\begin{dmath}
	\mathfrak H_N(\boldsymbol c) =
	c_1
	+
	\sum_{k=2}^{N}c_k^{1/k}\mathfrak u_k	
\end{dmath}
we obtain
\begin{dmath}\label{eq:universal-polynomial-linearization}
	\ee^{P_N(s; \boldsymbol c)}
	=
	\ee^{
	s
	\mathfrak H_N(\boldsymbol c)
	}
	\ground{\varphi(\boldsymbol t)}
\end{dmath}.

Given the exact correspondence~\eqref{eq:universal-equivalence}, the symbols
$\{\hat h_k, \varphi_0^{(k)}\}$ will not be used further, in favour of $\{\mathfrak u_k, \varphi_k(t_k)\}$, which enable a seamless transition to the analytic umbral framework.

\subsection{AUF realization}
\label{subsec:auf-brief}

The universal identity has an immediate analytic umbral translation. In the analytic umbral framework, the umbral operator $\mathfrak u$ is reinterpreted as a functional in the Borel additive complex variable u, namely:
\begin{dmath*}
	\mathfrak u \mapsto \ee^u	
\end{dmath*}.
For a suitable $F$ and $\zeta \in \mathbb C$, the formal operator function $F(\zeta\mathfrak u)$ thus corresponds to the composed \emph{umbral Borel functional} $F(\zeta\ee^u)$. After introducing the multiplicative Borel variable $z=\ee^u$, one defines the spectral \emph{jump function} associated to $F$ as
\begin{dmath*}
	J_F(t; \eta) \defeq {
		\mathcal M_z[F(\zeta z)](-t) =
		\int_0^\infty F(\zeta z)\,z^{-t-1}\,\mathrm d z
	}
\end{dmath*},
where $\mathcal M$ denotes the Mellin transform operator.

In particular, if $F(\zeta \ee^u) = \exp(\eta\zeta\ee^u)$, with $\eta \in \{-1, 1\}$, we have that
\begin{dmath*}
	J_{\mathrm{exp}_\eta}(t; \zeta) = (-\eta\zeta)^t \,\Gamma(-t)
\end{dmath*},
with the branch fixed by the principal determination of the logarithm. 

If $\varphi(t)$ is a meromorphic function with at most exponential growth on vertical lines, the analytic umbral pairing of an umbral Borel functional $F(\zeta \ee^u)$ with $\varphi(t)$ is defined by the Mellin--Barnes integral
\begin{dmath}\label{eq:umbral_analytic_pairing}
	f(\zeta) = {
		\langle F(\zeta \ee^u), \varphi \rangle 
		\defeq
		\frac{1}{2\pi i} \int_{\mathcal C} J_F(t; \zeta)\,\varphi(t)\,\mathrm d t
	}
\end{dmath},
where $\mathcal C$ is an upward oriented Barnes contour separating the poles of $J_F(t; \zeta)$ from the other singularities of the integrand.

In particular, the pairing with $\varphi$ of the exponential umbral Borel functional $\exp(\eta\zeta\ee^u)$ is
\begin{dmath*}
	\left\langle
		\exp(\eta\zeta\ee^u),
		\varphi(t)
	\right\rangle
	=
	\frac{1}{2\pi i}
	\int_{\mathcal C}
	(-\eta\zeta)^t
	\Gamma(-t)\,
	\varphi(t)
	\,\mathrm d t
\end{dmath*}.
We refer to \cite{RicciAUF2026} for further details. 

On an admissible Mellin--Barnes contour \(\mathcal C_k\), the AUF counterpart of \cref{eq:universal-exponential} is
\begin{dmath}\label{eq:single-auf-pairing}
\left\langle
\exp\!\left(
s\ee^{u_k}
\right),
\varphi_k(t_k)
\right\rangle
=
\frac{1}{2\pi\ii}
\int_{\mathcal C_k}
(-s)^{t_k}
\Gamma(-t_k)
\varphi_k(t_k)\,
\dd t_k
\end{dmath},
as can be verified by closing the contour to the right and taking the residues of \(\Gamma(-t_k)\). The polynomial linearization
therefore becomes the factorized multiple pairing
\begin{dmath}\label{eq:factorized-auf-pairing}
\exp\!\left(
P_N(s; \boldsymbol c)
\right)
=
\left\langle
\exp\!\left[
s\left(
c_1
+
\sum_{k=2}^{N}c_k^{1/k}\ee^{u_k}
\right)
\right],
\varphi(\boldsymbol t)
\right\rangle
\end{dmath}.
This is the AUF form of \cref{eq:universal-polynomial-linearization}. It does not, by
itself, authorize an interchange between a canonical \(s\)-integration and
the multiple Mellin--Barnes pairing; that issue is addressed for the Pearcey integral in
Section~\ref{sec:reversed}.

\section{Oscillatory canonical catastrophe integrals}
\label{sec:canonical}

For the \(A_k\) hierarchy, the normalized polynomial phase is
\begin{dmath}\label{eq:Ak-phase}
\Phi_{k+1}(s;\boldsymbol a)
= {
\frac{s^{k+1}}{k+1}
+
Q_{k-1}(s;\boldsymbol a),
\qquad
Q_{k-1}(s;\boldsymbol a)
=
\sum_{j=1}^{k-1}\frac{a_j}{j}s^j
}
\end{dmath}.
The coefficients
\(\boldsymbol a=(a_1,\ldots,a_{k-1})\) are the unfolding parameters.
The degree-\(k\) term is absent because it can be removed by a
translation of the integration variable.

The canonical field is conventionally written as the oscillatory
real-axis integral
\begin{dmath}\label{eq:Ak-field-formal}
\Psi_k(\boldsymbol a)
=
\int_{\mathbb R}^{\mathrm{osc}}
\exp\!\left[
\ii\Phi_{k+1}(s;\boldsymbol a)
\right]\dd s
\end{dmath}.
For real unfolding parameters, the two real half-line integrals
converge conditionally.  For the analysis below, it is more convenient
to use an equivalent contour on which the leading exponential decays.
This also defines the continuation to complex unfolding parameters.

Set \(m=k+1\).  Along \(s=r\ee^{\ii\theta}\),
\begin{dmath*}
\left|
\exp\!\left(\frac{\ii s^m}{m}\right)
\right|
=
\exp\!\left[
-\frac{r^m}{m}\sin(m\theta)
\right]
\end{dmath*}.
The decay sectors therefore satisfy \(\sin(m\theta)>0\), and their
central rays have arguments
\begin{dmath*}
\theta_j
= {
\frac{\pi/2+2\pi j}{m},
\qquad j\in\mathbb Z
}
\end{dmath*}.
On these rays, the leading exponential is simply
\(\exp(-r^m/m)\).

Let \(\Gamma_\theta\) denote the outward-oriented ray
\begin{dmath*}
\Gamma_\theta
= {
\{s=\ee^{\ii\theta}r:r\in[0,\infty)\}
}
\end{dmath*}
A contour entering from infinity along \(\theta_-\) and leaving
towards infinity along \(\theta_+\) is written
\begin{dmath*}
\Gamma_m=-\Gamma_{\theta_-}+\Gamma_{\theta_+}
\end{dmath*}.
The minus sign records the inward orientation of the first ray.

The standard contour is obtained by moving each end of the real axis
into its adjacent decay sector.  The positive end rotates into
\(0<\theta<\pi/m\).  At the negative end, the adjacent decay sector
lies above the real axis for odd \(m\), and below it for even \(m\).
Choosing the central rays gives
\begin{equation}\label{eq:prescribed-contour}
	\begin{aligned}
		\Gamma_m
		&=-\Gamma_{\pi-\pi/(2m)}+\Gamma_{\pi/(2m)},
		&\qquad &m=k+1\text{ odd},\\
		\Gamma_m
		&=-\Gamma_{\pi+\pi/(2m)}+\Gamma_{\pi/(2m)},
		&\qquad &m=k+1\text{ even}.
	\end{aligned}
\end{equation}
For even \(m\), this is the straight rotated line
\begin{dmath*}
	\Gamma_m=\ee^{\ii\pi/(2m)}\mathbb R
\end{dmath*}.
For odd \(m\), the two rays are not opposite.

We accordingly use the absolutely convergent contour representation
\begin{dmath}\label{eq:Ak-field}
\Psi_k(\boldsymbol a;\Gamma_{k+1})
=
\int_{\Gamma_{k+1}}
\exp\!\left[
\ii\Phi_{k+1}(s;\boldsymbol a)
\right]\dd s
\end{dmath}.
For real parameters it agrees with the standard oscillatory
real-axis integral; for complex parameters it defines its entire
continuation.  Moving the rays within the same decay sectors leaves
the field unchanged.  Choosing a different pair of endpoint sectors
generally selects a different solution of the canonical differential
system.

The Hermite representation follows by setting
\begin{dmath}\label{eq:Ak-roots}
c_j = {
\frac{\ii a_j}{j},
\qquad j=1,\ldots,k-1
}
\end{dmath}.
Then
\begin{dmath*}
P_{k-1}(s;\boldsymbol c)
=
\ii Q_{k-1}(s;\boldsymbol a)
\end{dmath*},
and the generating function gives
\begin{dmath}\label{eq:Ak-integrand-umbral}
	\begin{aligned}
	\exp\!\left[
	\ii\Phi_{k+1}(s;\boldsymbol a)
	\right]
	&=
	\exp\!\left(\frac{\ii s^{k+1}}{k+1}\right)
	\exp\!\left(P_{k-1}(s;\boldsymbol c)\right)
	\\
	&=
	\exp\!\left(\frac{\ii s^{k+1}}{k+1}\right)
	\exp\!\left(
	s\mathfrak H_{k-1}(\boldsymbol c)
	\right)
	\ground{\varphi(\boldsymbol t)}
	\\
	&=
	\exp\!\left(\frac{\ii s^{k+1}}{k+1}\right)
	\sum_{n=0}^{\infty}
	\frac{s^n}{n!}
	H_n^{[k-1]}(\boldsymbol c).
	\end{aligned}
\end{dmath}
Here the product ground state contains the components of orders
\(2,\ldots,k-1\).  For \(k=2\), the product is empty and
\(H_n^{[1]}(c_1)=c_1^n\).

Substitution into the contour integral yields
\begin{dmath}\label{eq:Ak-integral}
	\Psi_k(\boldsymbol a;\Gamma_m)
	= {
	\sum_{n=0}^{\infty}
	\frac{H_n^{[k-1]}(\boldsymbol c)}{n!}
	M_n^{(m)}(\Gamma_m),
	\qquad m=k+1
	}
\end{dmath},
where
\begin{dmath}\label{eq:moments}
M_n^{(m)}(\Gamma_m)
=
\int_{\Gamma_m}
s^n
\exp\!\left(\frac{\ii s^m}{m}\right)\dd s
\end{dmath}.

The interchange of sum and integral is justified directly by the
explicit coefficients of the Hermite polynomials.  They give the
bound
\begin{dmath}\label{eq:Hermite-absolute-majorant}
	\sum_{n=0}^{\infty}
	\frac{
	\left|H_n^{[m-2]}(\boldsymbol c)\right|
	}{n!}|s|^n
	\leq
	\exp\!\left(
	\sum_{j=1}^{m-2}|c_j|\,|s|^j
	\right)
\end{dmath}.
On the central rays in \cref{eq:prescribed-contour}, the modulus of
the leading exponential is \(\exp(-r^m/m)\).
For parameters in a compact set \(K\), the absolute series of
integrands is therefore bounded by
\begin{dmath*}
\exp\!\left[
-\frac{r^m}{m}+C_K(1+r^{m-2})
\right]
\end{dmath*},
which is integrable on \([0,\infty)\).  Hence
\cref{eq:Ak-integral} converges absolutely and locally uniformly in
the unfolding parameters.  In particular, it is a convergent
Hermite expansion, not merely an asymptotic series. Related Hermite
representations of diffraction integrals have also been developed in \cite{BelafhalEtAl2020,DattoliLicciardiSabia2025}.

The moments are elementary Gamma integrals.  Define, for
\(n\in\mathbb N_0\),
\begin{dmath}\label{eq:radial-canonical-moment}
J_n^{(m)}
= {
\int_0^\infty r^n\ee^{-r^m/m}\dd r
=
m^{\frac{n+1}{m}-1}
\Gamma\!\left(\frac{n+1}{m}\right)
}
\end{dmath}
On an outward central ray, the substitution
\(s=\ee^{\ii\theta_j}r\) gives
\begin{dmath}\label{eq:moments-on-a-ray}
\int_{\Gamma_{\theta_j}}
s^n\exp\!\left(\frac{\ii s^m}{m}\right)\dd s
=
\ee^{\ii(n+1)\theta_j}J_n^{(m)}
\end{dmath}.
The phase includes both \(s^n\) and the differential \(\dd s\).
Accounting for the orientation of the two rays, we obtain
\begin{dmath}\label{eq:two-ray-canonical-moments}
M_n^{(m)}(\Gamma_m)
=
J_n^{(m)}
\left[
\ee^{\ii(n+1)\theta_+}
-
\ee^{\ii(n+1)\theta_-}
\right].
\end{dmath}
Thus the radial Gamma factor is universal, while the bracket
records the selected contour.

For odd \(m\), substitution of the angles in
\cref{eq:prescribed-contour} gives
\begin{dmath}\label{eq:odd-degree-canonical-moments}
M_n^{(m)}(\Gamma_m)
=
J_n^{(m)}
\left[
\ee^{\frac{\ii\pi(n+1)}{2m}}
+
(-1)^n\ee^{-\frac{\ii\pi(n+1)}{2m}}
\right]
\end{dmath}.
This formula applies to the cubic \(A_2\) and quintic \(A_4\)
leading kernels.

For even \(m\), the same substitution gives
\begin{dmath}\label{eq:even-degree-oscillatory-moments}
M_n^{(m)}(\Gamma_m)
=
\left[1+(-1)^n\right]
\ee^{\frac{\ii\pi(n+1)}{2m}}
J_n^{(m)}
\end{dmath}.
Its relation to the Euclidean moment is particularly simple.
Since \(\Gamma_m=\ee^{\ii\pi/(2m)}\mathbb R\), put
\(s=\ee^{\ii\pi/(2m)}r\).  Then
\begin{dmath}\label{eq:oscillatory-Euclidean-moment}
M_n^{(m)}(\Gamma_m)
=
\ee^{\frac{\ii\pi(n+1)}{2m}}
M_{n,\mathrm E}^{(m)}
\end{dmath},
where
\begin{dmath}\label{eq:even-degree-Euclidean-moments}
M_{n,\mathrm E}^{(m)}
= {
\int_{\mathbb R}
r^n\ee^{-r^m/m}\dd r
=
\left[1+(-1)^n\right]J_n^{(m)}
}
\end{dmath}.
The factor \(1+(-1)^n\) expresses parity: odd moments vanish,
and even moments are twice the positive-half-line moment.

All moments used in \cref{eq:Ak-integral} are therefore defined
by absolutely convergent contour integrals.  Their real-axis
oscillatory notation requires a regularized interpretation for high
\(n\); no such regularization is needed on the prescribed decaying
contours.  The Hermite expansion and the Gamma-moment formulas above
are established directly on those contours
\cite[\href{https://dlmf.nist.gov/36}{(ch. 36)}]{NIST:DLMF}.

\section{Umbral representations of the canonical fields}
\label{sec:representations}

The Hermite expansion derived in the preceding section provides a
general construction of umbral representations for the integrated
catastrophe fields.  The essential step is to use the moments of the
leading kernel as the ground state of the universal umbral operator.
The polynomial structure of the unfolding is then retained either in
the powers of a single umbral operator or in a coupled multivariable ground
state.

\subsection{The canonical moment ground state}

Set $m=k+1$ and $d=k-1=m-2$.  Recall that
\begin{dmath}\label{eq:field-Hermite-moment-expansion}
	\Psi_k(\boldsymbol a;\Gamma_m)
	= {
	\sum_{n=0}^{\infty}
	\frac{H_n^{[d]}(\boldsymbol c)}{n!}
	M_n^{(m)}(\Gamma_m),
	\qquad
	c_j=\frac{\ii a_j}{j}
	}
\end{dmath}.
Here $M_n^{(m)}(\Gamma_m)$ is an ordinary contour moment, evaluated
after the leading umbral exponential has reconstructed the decaying
canonical kernel.

For a two-ray cycle
$\Gamma_m=-\Gamma_{\theta_-}+\Gamma_{\theta_+}$, introduce the
moment ground state
\begin{dmath}\label{eq:canonical-moment-ground-state}
	\mu_{m,\Gamma_m}(t)
	\defeq
	m^{\frac{t+1}{m}-1}
	\Gamma\!\left(\frac{t+1}{m}\right)
	\left[
	\ee^{\ii(t+1)\theta_+}
	-
	\ee^{\ii(t+1)\theta_-}
	\right]
\end{dmath}.
Its values on the nonnegative integers are precisely the canonical
moments:
\begin{dmath}\label{eq:moment-ground-state-lattice}
	\mu_{m,\Gamma_m}(n)
	= {
	M_n^{(m)}(\Gamma_m),
	\qquad n\in\mathbb N_0
	}
\end{dmath}.

The contour dependence of the field is therefore carried by the
moment ground state.  Different admissible cycles generally give
different functions $\mu_{m,\Gamma_m}$, although the radial
Gamma factor is the same.

At the formal level, only the lattice values in
\cref{eq:moment-ground-state-lattice} are needed.
Equation~\eqref{eq:canonical-moment-ground-state} supplies an explicit
meromorphic interpolation of those values.  Its use in an AUF pairing
requires the additional verification of Mellin admissibility and
vertical-growth conditions.

\subsection{Representation by a single universal umbra}

Let $\mathfrak u$ act on the canonical moment ground state according
to the universal evaluation rule
\begin{dmath}\label{eq:canonical-moment-evaluation}
	\mathfrak u^n
	\ground{\mu_{m,\Gamma_m}(t)}
	= {
	\mu_{m,\Gamma_m}(n),
	\qquad n\in\mathbb N_0
	}
\end{dmath}.
Substituting this rule into
\cref{eq:field-Hermite-moment-expansion} gives
\begin{dmath}\label{eq:Hermite-series-umbral-evaluation}
	\Psi_k(\boldsymbol a;\Gamma_m)
	=
	\sum_{n=0}^{\infty}
	\frac{H_n^{[d]}(\boldsymbol c)}{n!}
	\mathfrak u^n
	\ground{\mu_{m,\Gamma_m}(t)}
\end{dmath}.
The Hermite generating function now yields
\begin{dmath}\label{eq:canonical-single-umbral}
	\Psi_k(\boldsymbol a;\Gamma_{k+1})
	=
	\exp\!\left(
	\sum_{j=1}^{k-1}
	\frac{\ii a_j}{j}\mathfrak u^j
	\right)
	\ground{\mu_{k+1,\Gamma_{k+1}}(t)}
\end{dmath}.

In this representation, a single universal umbral operator carries the moments
of the leading catastrophe kernel.  Its powers
$\mathfrak u,\mathfrak u^2,\ldots,\mathfrak u^{k-1}$ retain the
degrees of the unfolding monomials.

The exponential in \cref{eq:canonical-single-umbral} is evaluated
through its generating series in $\mathfrak u$.  The resulting series
is exactly the convergent Hermite--moment expansion established in
the preceding section.  Thus this notation introduces no new
interchange of integration and umbral evaluation.

\subsection{Representation by a coupled multivariable ground state}

An equivalent representation makes the exponential linear in several
universal umbrae.  To obtain it, expand the Hermite coefficients:
\begin{dmath}\label{eq:Hermite-weighted-index-expansion}
	H_n^{[d]}(\boldsymbol c)
	=
	n!
	\sum_{\substack{
	r_1,\ldots,r_d\geq0\\
	r_1+2r_2+\cdots+d r_d=n}}
	\prod_{j=1}^{d}
	\frac{c_j^{r_j}}{r_j!}
\end{dmath}.
Substitution into
\cref{eq:field-Hermite-moment-expansion} and regrouping the absolutely
convergent series give
\begin{dmath}\label{eq:canonical-multiindex-series}
	\Psi_k(\boldsymbol a;\Gamma_m)
	=
	\sum_{r_1,\ldots,r_d\geq0}
	\left(
	\prod_{j=1}^{d}\frac{c_j^{r_j}}{r_j!}
	\right)
	\mu_{m,\Gamma_m}
	\left(
	r_1+2r_2+\cdots+d r_d
	\right)
\end{dmath}.

Define the coupled ground state
\begin{dmath}\label{eq:canonical-coupled-ground-state}
	\Phi_{m,\Gamma_m}(t_1,\ldots,t_d)
	\defeq
	\mu_{m,\Gamma_m}
	\left(
	t_1+2t_2+\cdots+d t_d
	\right).
\end{dmath}
Let $\mathfrak u_j$ denote the universal evaluation operator acting
in the variable $t_j$.  Their joint action is
\begin{dmath}\label{eq:canonical-joint-evaluation}
	\mathfrak u_1^{r_1}\cdots\mathfrak u_d^{r_d}
	\ground{\Phi_{m,\Gamma_m}(\boldsymbol t)}
	=
	\Phi_{m,\Gamma_m}(r_1,\ldots,r_d)
	=
	M_{r_1+2r_2+\cdots+d r_d}^{(m)}(\Gamma_m)
\end{dmath}.
Consequently,
\begin{dmath}\label{eq:canonical-multivariable-umbral}
	\Psi_k(\boldsymbol a;\Gamma_{k+1})
	=
	\exp\!\left(
	\sum_{j=1}^{k-1}
	\frac{\ii a_j}{j}\mathfrak u_j
	\right)
	\ground{\Phi_{k+1,\Gamma_{k+1}}(\boldsymbol t)}
\end{dmath}.

The operators $\mathfrak u_j$ commute, but the ground state does not
factorize.  Canonical integration couples the indices through the
weighted sum $t_1+2t_2+\cdots+d t_d$.  The equivalence with the
single-umbra representation is expressed by
\begin{dmath}\label{eq:single-multiple-evaluation-equivalence}
	\prod_{j=1}^{d}\mathfrak u_j^{r_j}
	\,
	\ground{\Phi_{m,\Gamma_m}(\boldsymbol t)}
	=
	\mathfrak u^{\sum_{j=1}^{d}j r_j}
	\ground{\mu_{m,\Gamma_m}(t)}
\end{dmath}.
This is an equality of evaluated expressions, not an identification
of the separate operators $\mathfrak u_j$ with powers of one
operator.

The two representations organize the same coefficient array in
different ways.  The single-operator form preserves the polynomial
unfolding in its exponent; the multivariable form transfers the
weighted polynomial structure into the ground state.

\subsection{The fold $A_2$}

For the normalized fold phase,
\begin{dmath*}
	\Phi_3(s;x)=\frac{s^3}{3}+xs
\end{dmath*},
the standard contour is
\begin{dmath*}
	\Gamma_3
	=
	-\Gamma_{5\pi/6}+\Gamma_{\pi/6}
\end{dmath*}.
Since the unfolding is linear, its Hermite polynomials reduce to
$H_n^{[1]}(\ii x)=(\ii x)^n$.  The general construction therefore gives
\begin{dmath}\label{eq:A2-moment-umbral}
	\Psi_2(x;\Gamma_3)
	= {
	\sum_{n=0}^{\infty}
	\frac{(\ii x)^n}{n!}
	\mu_{3,\Gamma_3}(n)
	=
	\ee^{\ii x\mathfrak u}
	\ground{\mu_{3,\Gamma_3}(t)}
	}
\end{dmath},
where the explicit moment ground state, given by \cref{eq:canonical-moment-ground-state}, can be written as
\begin{dmath}\label{eq:A2-moment-ground-state-simplified}
	\mu_{3,\Gamma_3}(t)
	=
	2\pi\,
	\ee^{\ii\pi t/2}
	\frac{
	3^{(t-2)/3}
	}{
	\Gamma((2-t)/3)
	}
\end{dmath}.
This function is entire: the apparent poles of the original
Gamma factor are cancelled by the contour phase difference.

To identify the field, first take $x$ real and return to the
equivalent oscillatory real-axis contour.  Since
$s^3/3+xs$ is odd in $s$, the two half-line contributions combine
into a cosine integral:
\begin{dmath}\label{eq:A2-Airy-identification}
	\Psi_2(x;\Gamma_3)
	= {
	2\int_0^\infty
	\cos\!\left(\frac{s^3}{3}+xs\right)\dd s
	=
	2\pi\operatorname{Ai}(x)
	}
\end{dmath}.
The last equality is the standard integral representation of
the Airy function \cite[\href{https://dlmf.nist.gov/9.5.E1}{eq. 9.5-E1}]{NIST:DLMF}.  Since both the contour-defined
field and $\operatorname{Ai}$ are entire, the identity extends
to every complex $x$.

The same identification is visible directly in the umbral series.
At integer $n$, the phase in
\cref{eq:A2-moment-ground-state-simplified} is
$\ee^{\ii\pi n/2}=\ii^n$.  Consequently,
\begin{dmath}\label{eq:A2-Airy-series}
	\operatorname{Ai}(x)
	=
	\sum_{n=0}^{\infty}
	\frac{(-x)^n}{n!}
	\frac{
	3^{(n-2)/3}
	}{
	\Gamma((2-n)/3)
	}
\end{dmath}.
In particular, the coefficients with $n=3q+2$ vanish, because
$1/\Gamma(-q)=0$ for $q\in\mathbb N_0$.  The remaining coefficients
form the two nonzero congruence classes of the Airy Taylor series.

If we wish to absorb the contour phase and the normalization
$2\pi$ into the ground state, define
\begin{dmath}\label{eq:Airy-normalized-ground-state}
	\varphi_{\mathrm{Ai}}(t)
	\defeq {
	\frac{
	\ee^{-\ii\pi t/2}
	}{
	2\pi
	}
	\mu_{3,\Gamma_3}(t)
	=
	\frac{
	3^{(t-2)/3}
	}{
	\Gamma((2-t)/3)
	}
	}
\end{dmath}.
The two equivalent formal representations are therefore
\begin{dmath}\label{eq:A2-Airy-umbral-equivalence}
	\boxed{
	2\pi\operatorname{Ai}(x)
	=
	\ee^{\ii x\mathfrak u}
	\llbracket\mu_{3,\Gamma_3}(t)\rrbracket,
	\qquad
	\operatorname{Ai}(x)
	=
	\ee^{-x\mathfrak u}
	\llbracket\varphi_{\mathrm{Ai}}(t)\rrbracket
	}
\end{dmath}.
The first retains the moments and phases of the canonical
oscillatory contour.  The second is normalized directly to the
Airy function and places the alternating sign in the exponential.

\subsection{The cusp $A_3$}

For the oscillatory Pearcey phase,
\begin{dmath*}
	\Phi_4(s;x,y)
	=
	{
	\frac{s^4}{4}+\frac{x}{2}s^2+ys,
	\qquad
	\Gamma_4=\ee^{\ii\pi/8}\mathbb R
	}
\end{dmath*}.
The corresponding field is the Pearcey integral
\cite{Pearcey1946},\cite[\href{https://dlmf.nist.gov/36.2.E14}{eq. 36.2.E14}]{NIST:DLMF}, in the normalization
\begin{dmath}\label{eq:Pearcey-definition}
	\operatorname{Pe}(x,y)
	=
	\int_{\Gamma_4}
	\exp\!\left[
	\ii\left(
	\frac{s^4}{4}+\frac{x}{2}s^2+ys
	\right)
	\right]\dd s
\end{dmath}.
If $P(X,Y)$ denotes the convention
$\int_{\mathbb R}^{\mathrm{osc}}
\exp[\ii(t^4+Xt^2+Yt)]\dd t$
used in \cite[Eq.~36.2.14]{NIST:DLMF}, the change
$s=\sqrt{2}\,t$ gives
\begin{dmath*}
	\operatorname{Pe}(x,y)
	=
	\sqrt{2}\,P(x,\sqrt{2}\,y)
\end{dmath*}.

The moment ground state is
\begin{dmath}\label{eq:A3-moment-ground-state}
	\mu_{4,\Gamma_4}(t)
	=
	4^{\frac{t+1}{4}-1}
	\Gamma\!\left(\frac{t+1}{4}\right)
	\left[
	\ee^{\ii\pi(t+1)/8}
	-
	\ee^{9\ii\pi(t+1)/8}
	\right]
\end{dmath},
or, equivalently,
\begin{dmath}\label{eq:A3-moment-ground-state-parity}
	\mu_{4,\Gamma_4}(t)
	=
	4^{\frac{t+1}{4}-1}
	\Gamma\!\left(\frac{t+1}{4}\right)
	\ee^{\ii\pi(t+1)/8}
	\left(1+\ee^{\ii\pi t}\right)
\end{dmath},
with removable singularities filled by continuation.
At nonnegative integers, this gives
\begin{dmath*}
	{
	\mu_{4,\Gamma_4}(2n)
	=
	2^{n-\frac12}
	\Gamma\!\left(\frac{2n+1}{4}\right)
	\ee^{\ii\pi(2n+1)/8},
	\qquad
	\mu_{4,\Gamma_4}(2n+1)=0
	}
\end{dmath*}.

The two general umbral representations become
\begin{dmath}\label{eq:A3-two-umbral-forms}
	\operatorname{Pe}(x,y)
	=
	\exp\!\left(
	\ii y\mathfrak u+\frac{\ii x}{2}\mathfrak u^2
	\right)
	\ground{\mu_{4,\Gamma_4}(t)}
	=
	\exp\!\left(
	\ii y\mathfrak u_1+\frac{\ii x}{2}\mathfrak u_2
	\right)
	\ground{\Phi_{4,\Gamma_4}(t_1,t_2)}
\end{dmath},
where
\begin{dmath*}
	\Phi_{4,\Gamma_4}(t_1,t_2)
	=
	\mu_{4,\Gamma_4}(t_1+2t_2)
\end{dmath*}.
The first index in this coupled ground state counts powers of
$y$, while the second counts powers of $x$.  The vanishing odd
moments imply
$\operatorname{Pe}(x,-y)=\operatorname{Pe}(x,y)$.

To explain the Euclidean representation, first retain the linear
source variable and define
\begin{dmath}\label{eq:Euclidean-Pearcey-linear-source}
	\operatorname{Pe}_E(x,y)
	\defeq
	\int_{\mathbb R}
	\exp\!\left(
	-\frac{s^4}{4}-\frac{x}{2}s^2-ys
	\right)\dd s
	=
	\int_{\mathbb R}
	\ee^{-s^4/4-xs^2/2}\cosh(ys)\dd s
\end{dmath}.
The second equality follows by symmetry under $s\mapsto-s$.
In particular, $\operatorname{Pe}_E$ is even in $y$.

The rotation $s=\ee^{\ii\pi/8}r$ in
\cref{eq:Pearcey-definition} gives
\begin{dmath}\label{eq:Pearcey-Euclidean-rotation}
	\boxed{
	\operatorname{Pe}(x,y)
	=
	\ee^{\ii\pi/8}
	\operatorname{Pe}_E
	\left(
	\ee^{-\ii\pi/4}x,
	\ee^{-3\ii\pi/8}y
	\right)
	}
\end{dmath}.
The prefactor is the Jacobian of the rotation; the rotated
arguments follow from the quadratic and linear terms of the phase.

The Euclidean integral is absolutely convergent for every
$(x,y)\in\mathbb C^2$.  On compact parameter sets its integrand
is bounded by
\begin{dmath*}
	\exp\!\left(
	-\frac{s^4}{4}+C_Ks^2+C_K|s|
	\right)
\end{dmath*},
and the same bound, multiplied by a polynomial in $s$, applies
after parameter differentiation.  Hence
$\operatorname{Pe}_E(x,y)$ is jointly entire in $x,y$.
Equation~\eqref{eq:Pearcey-Euclidean-rotation} establishes the
same property for the contour-defined oscillatory field.
For complex parameters, it is this contour representation,
rather than an everywhere-convergent real-axis oscillatory
integral, that defines $\operatorname{Pe}$.

Although \cref{eq:A3-two-umbral-forms} already gives complete
formal representations, the evenness in the source suggests a
more economical indexing.  Introduce $\eta=y^2$ and define
\begin{dmath}\label{eq:PE-definition}
	\mathcal P_E(x,\eta)
	\defeq {
	\int_{\mathbb R}
	\ee^{-s^4/4-xs^2/2}
	\cosh(s\sqrt{\eta})\dd s,
	\qquad
	\operatorname{Pe}_E(x,y)=\mathcal P_E(x,y^2)
	}
\end{dmath}.
There is no square-root ambiguity: $\cosh(s\sqrt{\eta})$ is
entire in $\eta$, as its power series contains only integer
powers of $\eta$.

Expanding the quadratic and source factors yields
\begin{dmath*}
	\mathcal P_E(x,\eta)
	=
	\sum_{p=0}^{\infty}\sum_{q=0}^{\infty}
	\frac{(-x/2)^p\eta^q}{p!(2q)!}
	\int_{\mathbb R}
	s^{2p+2q}\ee^{-s^4/4}\dd s
\end{dmath*}.
Using
\begin{dmath*}
	\int_{\mathbb R}
	s^{2p+2q}\ee^{-s^4/4}\dd s
	=
	2^{p+q-\frac12}
	\Gamma\!\left(\frac{2p+2q+1}{4}\right)
\end{dmath*}
and the Gamma duplication formula, we obtain
\begin{dmath}\label{eq:PE-integrated-series}
	\mathcal P_E(x,\eta)
	=
	\sum_{p=0}^{\infty}\sum_{q=0}^{\infty}
	\frac{(-x)^p\eta^q}{p!\,q!}
	\Phi(p,q)
\end{dmath},
with
\begin{dmath}\label{eq:PE-ground-state}
	\Phi(t_1,t_2)
	=
	\sqrt{\frac{\pi}{2}}\,
	2^{-t_2}
	\frac{
	\Gamma\!\left((2t_1+2t_2+1)/4\right)
	}{
	\Gamma(t_2+1/2)
	}
\end{dmath}.
Hence the parity-reduced umbral representation is
\begin{dmath}\label{eq:PE-integrated-umbral}
	\mathcal P_E(x,\eta)
	=
	\exp\!\left(
	-x\mathfrak u_1+\eta\mathfrak u_2
	\right)
	\ground{\Phi(t_1,t_2)}
\end{dmath}.

The ground state $\Phi$ is introduced to encode only the
nonvanishing source coefficients.  Its first index counts
powers of $x$ and its second counts powers of $\eta=y^2$.
It therefore reorganizes the same quartic moments, including
the change from $(2q)!$ to $q!$ in the source coefficients.
Unlike $\Phi_{4,\Gamma_4}$, it does not need an explicit
evenness factor on the nonnegative integer lattice.

Finally, combining the parity reduction with the contour
rotation gives
\begin{dmath}\label{eq:oscillatory-Pearcey-umbral}
	\operatorname{Pe}(x,y)
	=
	\ee^{\ii\pi/8}
	\mathcal P_E
	\left(
	\ee^{-\ii\pi/4}x,
	\ee^{-3\ii\pi/4}y^2
	\right)
	=
	\ee^{\ii\pi/8}
	\exp\!\left(
	-\ee^{-\ii\pi/4}x\mathfrak u_1
	+
	\ee^{-3\ii\pi/4}y^2\mathfrak u_2
	\right)
	\ground{\Phi(t_1,t_2)}
\end{dmath}.

The ground state $\Phi(t,r)$ is also used in
the bivariate Mellin--Barnes AUF representation of $\mathcal P_E(x,\eta)$.
It is holomorphic in the chamber
\begin{dmath}\label{eq:PE-admissible-chamber}
	\RE t< {0,
	\qquad
	\RE r<0,
	\qquad
	\RE(t+r)>-\frac12
	}
\end{dmath}.
Indeed, the argument of its numerator Gamma function has
positive real part there, while the reciprocal Gamma function
in the denominator introduces no poles.

To check convergence of the AUF pairing, set
$t=c+\ii u$ and $r=d+\ii v$, with $c,d$ satisfying
\cref{eq:PE-admissible-chamber}. Stirling's estimate on vertical
lines gives
\begin{dmath}\label{eq:PE-vertical-bound}
	\left|
	\Gamma(-t)\Gamma(-r)\Phi(t,r)
	\right|
	\leq
	C(1+|u|+|v|)^N
	\exp\!\left(
	-\frac{\pi}{2}|u|
	-\frac{\pi}{4}|u+v|
	\right)
\end{dmath},
for suitable constants $C>0$ and $N\geq0$.
Here the exponential factors associated with
$\Gamma(-r)$ and $1/\Gamma(r+1/2)$ cancel, whereas the
coupled numerator contributes
$\exp(-\pi|u+v|/4)$. Since
$|u|+|u+v|$ controls $|u|+|v|$, the remaining decay is
sufficient for absolute integration over both vertical lines.

In particular, on the Fourier slice $\eta<0$, the AUF pairing
\begin{dmath}\label{eq:PE-admissible-pairing}
	\croc{\exp(-x\ee^{u_1} + \eta\ee^{u_2}), \Phi(t,r)} =  
	\frac{1}{(2\pi\ii)^2}
	\int_{c-\ii\infty}^{c+\ii\infty}
	\int_{d-\ii\infty}^{d+\ii\infty}
	x^t (-\eta)^r
	\Gamma(-t)\Gamma(-r)\Phi(t,r)
	\dd r\,\dd t
\end{dmath}
is absolutely convergent.

The convergence statement concerns the fixed vertical-contour
integral, not the full domain of the represented function.
The integral converges absolutely and locally
uniformly for nonzero $x$ and $\eta$ in the sector
\begin{dmath*}
	|\arg(-\eta)|<\frac{\pi}{4},
	\qquad
	|\arg x-\arg(-\eta)|<\frac{\pi}{2}
\end{dmath*},
on the selected branches. In particular, it converges for
$x>0$ and $\eta<0$. Positive $\eta$ lies outside this sector,
and zero arguments must be treated by continuation or limits.

To identify the function defined by the pairing, one first
evaluates \cref{eq:PE-vertical-bound} for $x>0$ and $\eta<0$, obtaining
$\mathcal P_E(x,\eta)$. Independently, the Euclidean integral
defines a jointly entire function. Holomorphic uniqueness therefore shows that the AUF pairing
represents the restriction of this entire function to its
convergence sector. Its value consequently has a unique
jointly entire continuation to all $(x,\eta)\in\mathbb C^2$,
although the original fixed-contour Mellin--Barnes integral
need not converge outside that sector.

\subsubsection{Quartic and Higgs-type partition functions}

The Euclidean Pearcey field also describes the partition function
of a zero-dimensional quartic scalar model with source term,
\begin{dmath}\label{eq:quartic-partition}
	Z_\lambda(a,b)
	= {
	\int_{\mathbb R}
	\exp\!\left(
	-\frac{\lambda}{4}\phi^4
	-\frac{a}{2}\phi^2
	-b\phi
	\right)\dd\phi,
	\qquad \lambda>0
	}
\end{dmath}.
For fixed $\lambda>0$, quartic decay makes this integral absolutely
convergent and jointly entire in $a$ and $b$. Its dependence on $b$
is even, as follows from the substitution $\phi\mapsto-\phi$.

The rescaling $s=\lambda^{1/4}\phi$ gives
\begin{dmath}\label{eq:quartic-pearcey-reduction}
	Z_\lambda(a,b)
	=
	\lambda^{-1/4}
	\mathcal P_E\!\left(
	\frac{a}{\sqrt{\lambda}},
	\frac{b^2}{\sqrt{\lambda}}
	\right)
\end{dmath}.
Consequently, no additional ground state is required: the
parity-reduced Pearcey ground state in
\cref{eq:PE-ground-state} gives directly
\begin{dmath}\label{eq:quartic-partition-umbral}
	Z_\lambda(a,b)
	=
	\lambda^{-1/4}
	\exp\!\left(
	-\frac{a}{\sqrt{\lambda}}\mathfrak u_1
	+\frac{b^2}{\sqrt{\lambda}}\mathfrak u_2
	\right)
	\ground{\Phi(t_1,t_2)}
\end{dmath}.
The corresponding AUF representation is obtained by the same
rescaling in \cref{eq:PE-admissible-pairing}, initially in its domain
of convergence and elsewhere by analytic continuation.

For real $a$ and $b$, the geometry of the model is determined by
the potential
\begin{dmath}\label{eq:higgs-potential}
	V(\phi;a,b)
	=
	\frac{\lambda}{4}\phi^4
	+\frac{a}{2}\phi^2
	+b\phi
\end{dmath}.
Its stationary points satisfy
\begin{dmath*}
	\lambda\phi^3+a\phi+b=0
\end{dmath*}.
Requiring two stationary points to coalesce, that is,
$V'(\phi)=V''(\phi)=0$, and eliminating $\phi$ gives the cusp
discriminant
\begin{dmath}\label{eq:higgs-discriminant}
	4a^3+27\lambda b^2=0
\end{dmath}.
This curve separates the region with one real stationary point
from that with three; it is a bifurcation set of the potential,
not a singularity of the exact partition function \cite[\href{https://dlmf.nist.gov/36.4}{sec. 36.4}]{NIST:DLMF}.

On the symmetric section $b=0$, the origin is the unique minimum
for $a>0$. For $a<0$, it becomes a local maximum and two minima
appear at
\begin{dmath*}
	\phi_\pm = {
	\pm v,
	\qquad
	v=\sqrt{-a/\lambda}
	}
\end{dmath*}.
In this regime the potential can be written as
\begin{dmath*}
	V(\phi;a,0)
	=
	\frac{\lambda}{4}(\phi^2-v^2)^2
	-\frac{a^2}{4\lambda}
\end{dmath*}.
This double-well structure motivates the designation
\emph{Higgs-type}; the source $b$ tilts the potential and removes
the degeneracy of its minima.

The analogy concerns the quartic potential and its stationary-point
geometry. It does not imply spontaneous symmetry breaking in the
finite zero-dimensional model: at zero source the exact expectation
value remains
\begin{dmath*}
	\langle\phi\rangle_{a,0}
	= {
	-\left.
	\frac{\partial}{\partial b}\log Z_\lambda(a,b)
	\right|_{b=0}
	=
	0
	}
\end{dmath*}.
Nor does this single real-variable integral include the field,
gauge or measure structure of the physical Higgs theory.
Its role here is to show that quartic scalar partition functions
and cusp diffraction fields share the same integrated umbral
ground state, although their physical interpretations differ.

\subsection{The swallowtail $A_4$}

For the normalized swallowtail phase,
\begin{dmath*}
	\Phi_5(s;a,b,c)
	=
	\frac{s^5}{5}
	+\frac{a}{3}s^3
	+\frac{b}{2}s^2
	+cs
\end{dmath*},
the standard oscillatory cycle is represented by
\begin{dmath*}
	\Gamma_5
	=
	-\Gamma_{9\pi/10}+\Gamma_{\pi/10}
\end{dmath*}.
Both rays provide quintic decay. The resulting field
$\mathcal S(a,b,c)$ is jointly entire in its three unfolding
parameters and agrees with the real-axis oscillatory integral
when these parameters are real.

The canonical moments are interpolated by
\begin{dmath}\label{eq:A4-moment-ground-state}
	\mu_{5,\Gamma_5}(t)
	=
	5^{(t-4)/5}
	\Gamma\!\left(\frac{t+1}{5}\right)
	\left[
	\ee^{\ii\pi(t+1)/10}
	-
	\ee^{9\ii\pi(t+1)/10}
	\right]
\end{dmath},
with removable singularities filled in. In particular,
$\mu_{5,\Gamma_5}(n)$ is the $n$th canonical moment for every
nonnegative integer $n$. The phase difference records the
orientations and directions of the two rays.

The single-umbra representation is therefore
\begin{dmath}\label{eq:A4-single-umbral}
	\mathcal S(a,b,c)
	=
	\exp\!\left(
	\ii c\mathfrak u
	+\frac{\ii b}{2}\mathfrak u^2
	+\frac{\ii a}{3}\mathfrak u^3
	\right)
	\ground{\mu_{5,\Gamma_5}(t)}
\end{dmath}.
Equivalently, introducing one umbral variable for each unfolding
parameter gives
\begin{dmath}\label{eq:A4-multivariable-umbral}
	\mathcal S(a,b,c)
	=
	\exp\!\left(
	\ii c\mathfrak u_1
	+\frac{\ii b}{2}\mathfrak u_2
	+\frac{\ii a}{3}\mathfrak u_3
	\right)
	\ground{\Phi_{5,\Gamma_5}(t_1,t_2,t_3)}
\end{dmath},
where
\begin{dmath}\label{eq:A4-coupled-ground-state}
	\Phi_{5,\Gamma_5}(t_1,t_2,t_3)
	=
	\mu_{5,\Gamma_5}(t_1+2t_2+3t_3)
\end{dmath}.
The weighted combination of indices reflects the linear,
quadratic and cubic unfolding monomials. These formulas are
understood through their convergent evaluated Taylor series;
an AUF realization additionally requires admissible Mellin
contours and growth estimates.

Unlike the quartic case, the quintic exponential
$\ee^{-s^5/5}$ does not decay at both ends of the real axis.
Thus a Euclidean realization must retain two distinct contour
rays rather than a single real-line integral.

\subsection{Higher catastrophes and analytic completion}

Equations~\eqref{eq:canonical-single-umbral} and
\eqref{eq:canonical-multivariable-umbral} apply throughout the
$A_k$ hierarchy.  Once the contour moments are known, no additional
formal construction is needed for higher $k$.  The increasing
number of unfolding parameters enters through the Hermite
polynomials and the weighted argument of the coupled ground state.

This construction follows the natural order throughout: the
leading kernel is evaluated, its ordinary contour moments are
computed, and the resulting coefficient array is expressed using
the universal umbral operator.  It does not require an integral
of an unevaluated umbral exponential.

The passage to an AUF representation is a further analytic step.
The lattice values determine the formal field, but do not uniquely
determine its interpolation away from the integer lattice.
An admissible Mellin representation must also specify branches,
contours and growth conditions.  In particular, an interpolation
that is convenient for formal evaluation need not itself provide
a convergent Mellin--Barnes pairing.  The Euclidean Pearcey ground
state in \cref{eq:PE-ground-state} provides the principal explicit
example developed below.

\section{Reversing evaluation and integration}
\label{sec:reversed}

The preceding constructions follow the natural order: evaluate the
umbral exponential and then integrate the resulting polynomial
exponential. We now ask whether integration can precede evaluation.
For the Euclidean Pearcey field, the AUF supplies an analytic meaning
to this reversed order that is not provided by the formal umbral framework alone.

\subsection{The limitation of formal evaluation}

Work initially with $x,z>0$ and set
\begin{dmath}\label{eq:G-definition}
	G(x,z)
	= {
	\mathcal P_E(x,-z)
	=
	2\int_0^\infty
	\ee^{-s^4/4-xs^2/2}
	\cos(s\sqrt z)\dd s
	}
\end{dmath}.
Choose
\begin{dmath*}
	a_2 = {
	\ii\sqrt{\frac{x}{2}},
	\qquad
	a_4=2^{-1/2}\ee^{3\pi\ii/4},
	\qquad
	\mathfrak A_x=a_2\mathfrak u_2+a_4\mathfrak u_4
	}
\end{dmath*}.
The formal representation is
\begin{dmath}\label{eq:formal-Pearcey-kernel}
	\ee^{-s\mathfrak A_x}
	\ground{\varphi_2(t_2)\varphi_4(t_4)}
	=
	\ee^{-xs^2/2-s^4/4}
\end{dmath}.
Thus evaluation before integration gives
\cref{eq:G-definition} without regularization.

For a complex number $A$, the elementary identity
\begin{dmath}\label{eq:scalar-resolvent}
	2\int_0^\infty
	\ee^{-sA}\cos(s\sqrt z)\dd s
	= {
	\frac{2A}{A^2+z},
	\qquad \RE A>0
	}
\end{dmath}
suggests replacing $A$ by $\mathfrak A_x$. However, the condition
$\RE\mathfrak A_x>0$ has no meaning in the formal algebra, and
no operator-valued integration has been defined there.
Adjoining the inverse of $\mathfrak A_x^2+z$ defines the rational
symbol
\begin{dmath*}
	\mathfrak R_z
	=
	\frac{2\mathfrak A_x}{\mathfrak A_x^2+z}
\end{dmath*},
but does not prove that it represents the integral of the
unevaluated exponential, or specify its action on the ground state.

Indeed, expanding the rational symbol in nonnegative powers of
$\mathfrak A_x$ gives
\begin{dmath}\label{eq:positive-power-resolvent}
	\mathfrak R_z
	=
	2\sum_{n=0}^{\infty}
	\frac{(-1)^n}{z^{n+1}}\mathfrak A_x^{2n+1}
\end{dmath}.
However,
\begin{dmath}\label{eq:odd-composite-moments}
	\mathfrak A_x^{2n+1}
	\ground{\varphi_2(t_2)\varphi_4(t_4)}
	=
	0
\end{dmath},
because a surviving monomial must contain an even power of
$\mathfrak u_2$ and a power of $\mathfrak u_4$ divisible by four.
Coefficientwise evaluation of this expansion therefore gives
zero, rather than the Pearcey field. This does not rule out
every extension of formal evaluation to rational expressions;
it shows that the positive-power prescription does not
provide the required extension.

The alternative expansion
\begin{dmath*}
	\mathfrak R_z
	=
	2\sum_{n=0}^{\infty}
	(-z)^n\mathfrak A_x^{-2n-1}
\end{dmath*}
requires an additional evaluation rule for negative powers of
the \emph{combined} umbral operator $\mathfrak A_x$. Definitions of negative powers of
the separate umbral operators $\mathfrak u_2$ and $\mathfrak u_4$ do not by themselves supply this rule.
Likewise, a formal Laplace transform produces a regulated resolvent
as a series at infinite regulator, but does not justify removing
the regulator in an unevaluated operator integral.

The obstruction is therefore not to the correctness of the formal umbral framework
in its natural order. Reversing the order requires an extension of the formal evaluation supplied by an analytic realization of the integrated kernel. Without such an analytic extension, substitution of \(\mathfrak A_x\) into \cref{eq:scalar-resolvent} remains an operational prescription rather than an established integral identity.

\subsection{The AUF extension and the Pearcey identity}

Our objective is now the evaluation of
$G(x,z)=\mathcal P_E(x,-z)$ in the analytic framework. Its integrand has the AUF
representation
\begin{dmath}\label{eq:auf-Pearcey-kernel}
	\ee^{-s^4/4-xs^2/2}
	=
	\left\langle
	\ee^{-s\mathcal A_x(\boldsymbol u)},
	\varphi_2(t_2)\varphi_4(t_4)
	\right\rangle,
	\qquad
	\mathcal A_x(\boldsymbol u)
	=
	a_2\ee^{u_2}+a_4\ee^{u_4}
\end{dmath}.
Explicitly
\begin{dmath}\label{eq:explicit-factorized-AUF-pairing}
	\left\langle
	\ee^{-s\mathcal A_x},
	\varphi_2\varphi_4
	\right\rangle
	=
	\frac{1}{(2\pi\ii)^2}
	\int_{\mathcal C_{t_2}}
	(sa_2)^{t_2}\Gamma(-t_2)\varphi_2(t_2)\dd t_2
	\int_{\mathcal C_{t_4}}
	(sa_4)^{t_4}\Gamma(-t_4)\varphi_4(t_4)\dd t_4
\end{dmath}.
The monomial ground states used in this calculation are
\begin{dmath}\label{eq:quadratic-quartic-ground-states}
	\varphi_2(\tau)
	= {
	\frac12\ee^{-\pi\ii\tau/2}
	\frac{\Gamma(-\tau/2)}{\Gamma(-\tau)},
	\qquad
	\varphi_4(\tau)
	=
	\frac14\ee^{-3\pi\ii\tau/4}
	\frac{\Gamma(-\tau/4)}{\Gamma(-\tau)}
	}
\end{dmath},
with removable singularities filled by continuation.
These are equivalent forms of the ground states introduced
earlier through the congruence selectors. Their explicit Gamma-quotient form also explains their
analytic role. In the AUF pairing, multiplication by
$\varphi_k(\tau)$ replaces the exponential Mellin density
$\Gamma(-\tau)$ by a rescaled Gamma factor:
\begin{dmath*}
	\Gamma(-\tau)\varphi_k(\tau)
	= {
	\frac1k
	\ee^{-\pi\ii(k-1)\tau/k}
	\Gamma\!\left(-\frac{\tau}{k}\right),
	\qquad k=2,4
	}
\end{dmath*}.
The choices
$a_2=\ii\sqrt{x/2}$ and
$a_4=2^{-1/2}\ee^{3\pi\ii/4}$
cancel the displayed phases. The substitutions
$v=\tau/2$ and $v=\tau/4$ then reduce the respective
pairings to the ordinary inverse Mellin formulas for
$\ee^{-xs^2/2}$ and $\ee^{-s^4/4}$. Other root choices are equivalent for formal evaluation,
but their Mellin realizations require corresponding
branch and admissibility checks.

The selected branches therefore turn each monomial pairing
into an ordinary inverse Mellin representation of a decaying
exponential. Their product reconstructs the Euclidean
Pearcey kernel. This is the concrete analytic meaning of
the formal monomial evaluation in the present calculation.

Evaluation before integration gives
\begin{dmath*}
	G(x,z)
	=
	2\int_0^\infty
	\left\langle
	\ee^{-s\mathcal A_x(\boldsymbol u)},
	\varphi_2(t_2)\varphi_4(t_4)
	\right\rangle
	\cos(s\sqrt z)\dd s
\end{dmath*}.

To investigate the reversed order, choose paths for the complex
logarithmic variables $u_2,u_4$ such that
$y_k=a_k\ee^{u_k}>0$. For $x>0$ and the selected branches,
this is achieved by
\begin{dmath*}
	u_2 = {
	v_2-\frac{\pi\ii}{2},
	\qquad
	u_4=v_4-\frac{3\pi\ii}{4},
	\qquad v_2,v_4\in\mathbb R
	}
\end{dmath*},
so that
\begin{dmath*}
	y_2 = {
	a_2\ee^{u_2}
	=
	\sqrt{\frac{x}{2}}\,\ee^{v_2},
	\qquad
	y_4=a_4\ee^{u_4}
	=
	2^{-1/2}\ee^{v_4}
	}
\end{dmath*}.
Both variables range independently over the positive real
axis. For $z>0$, integration of the transmutation kernel
is therefore an ordinary absolutely convergent integral:
\begin{dmath}\label{eq:integrated-transmutation-kernel}
	K_z(y_2,y_4)
	= {
	2\int_0^\infty
	\ee^{-s(y_2+y_4)}\cos(s\sqrt z)\dd s
	=
	\frac{2(y_2+y_4)}{(y_2+y_4)^2+z},
	\qquad y_2,y_4>0
	}
\end{dmath}.

This rational kernel is not itself the Pearcey field:
the ground-state pairing must still act on it. The identity
to be established is
\begin{dmath}\label{eq:kernel-Pearcey-target}
	G(x,z)
	\overset{\mathrm{?}}{=}
	\left\langle
	K_z\!\left(a_2\ee^{u_2},a_4\ee^{u_4}\right),
	\varphi_2(t_2)\varphi_4(t_4)
	\right\rangle
\end{dmath}.

The chosen paths justify the scalar integration defining
$K_z$; they do not by themselves justify exchanging
integration and evaluation. To establish that exchange,
we compute the Mellin transform of $K_z$ on the positive
quadrant and use it to define the right-hand pairing,
retaining the branch factors $a_2^{\tau_2}a_4^{\tau_4}$.
The following calculation proves that this pairing converges
and reproduces the independently obtained Mellin--Barnes
representation of $G(x,z)$.

The Mellin transform follows by setting
$A=y_2+y_4$ and $v=y_2/A$. The $v$-integration gives a beta
function, and the remaining $A$-integration yields
\begin{dmath}\label{eq:rational-kernel-mellin}
	\mathcal M [K_z](\sigma_2,\sigma_4)
	= {
	\frac{\Gamma(\sigma_2)\Gamma(\sigma_4)}
	{\Gamma(\Sigma)}
	z^{(\Sigma-1)/2}
	\Gamma\!\left(\frac{\Sigma+1}{2}\right)
	\Gamma\!\left(\frac{1-\Sigma}{2}\right),
	\qquad
	\Sigma=\sigma_2+\sigma_4
	}
\end{dmath},
in the nonempty chamber
\begin{dmath}\label{eq:rational-kernel-chamber}
	\RE\sigma_2 > {
	0,
	\qquad
	\RE\sigma_4>0,
	\qquad
	\RE\Sigma<1
	}
\end{dmath}.

To apply the AUF prescription, the Mellin transform must be
taken in the exponential variables $w_j=\ee^{u_j}$.
Define
\begin{dmath*}
	H_{x,z}(w_2,w_4)
	=
	K_z(a_2w_2,a_4w_4)
\end{dmath*}.
The jump function entering the pairing is therefore the
Mellin transform of $H_{x,z}$, evaluated at the negative
spectral indices, rather than the unscaled Mellin transform
of $K_z$.

For the complex coefficients chosen above, use the rays
\begin{dmath*}
	\ell_2 = {
	a_2^{-1}\mathbb R_+,
	\qquad
	\ell_4=a_4^{-1}\mathbb R_+
	}
\end{dmath*},
oriented from zero to infinity. These are precisely the
paths for which $y_j=a_jw_j>0$. Choose logarithms consistently
so that $\log w_j=\log y_j-\log a_j$ on $\ell_j$.
The corresponding sectorial Mellin transform is
\begin{dmath}\label{eq:scaled-kernel-mellin}
	\mathcal M_{\boldsymbol\ell}[H_{x,z}]
	(\sigma_2,\sigma_4)
	=
	\int_{\ell_2}\int_{\ell_4}
	w_2^{\sigma_2-1}w_4^{\sigma_4-1}
	K_z(a_2w_2,a_4w_4)
	\dd w_4\,\dd w_2
\end{dmath}.
Under $y_j=a_jw_j$, the powers and differentials satisfy
\begin{dmath*}
	w_j^{\sigma_j-1}\dd w_j
	=
	a_j^{-\sigma_j}
	y_j^{\sigma_j-1}\dd y_j
\end{dmath*}.
Consequently,
\begin{dmath}\label{eq:Mellin-scaling-identity}
	\mathcal M_{\boldsymbol\ell}[H_{x,z}]
	(\sigma_2,\sigma_4)
	=
	a_2^{-\sigma_2}a_4^{-\sigma_4}
	\mathcal M[K_z](\sigma_2,\sigma_4)
\end{dmath},
with the same convergence chamber
\cref{eq:rational-kernel-chamber}.
This is a direct change of variables on the selected rays;
no deformation across the poles of the rational kernel
is being assumed.

The AUF jump function is thus
\begin{dmath}\label{eq:resolvent-jump-function}
	J_{x,z}(\tau_2,\tau_4)
	=
	\mathcal M_{\boldsymbol\ell}[H_{x,z}]
	(-\tau_2,-\tau_4)
	=
	a_2^{\tau_2}a_4^{\tau_4}
	\mathcal M[K_z](-\tau_2,-\tau_4)
\end{dmath}.
Pairing this function with the product ground state gives
\begin{dmath}\label{eq:factorized-auf-resolvent}
	\mathcal R(x,z)
	=
	\frac{1}{(2\pi\ii)^2}
	\int_{\mathcal C_{\tau_2}}
	\int_{\mathcal C_{\tau_4}}
	J_{x,z}(\tau_2,\tau_4)
	\varphi_2(\tau_2)\varphi_4(\tau_4)
	\dd\tau_4\,\dd\tau_2
\end{dmath},
where
\begin{dmath*}
	\RE\tau_2< {
	0,
	\qquad
	\RE\tau_4<0,
	\qquad
	-1<\RE(\tau_2+\tau_4)<0
	}
\end{dmath*}.
This is the AUF prescription applied to the integrated
kernel, with its logarithmic paths and branches made
explicit. The remaining calculation establishes convergence
of this pairing and identifies its value with $G(x,z)$.
The explicit monomial ground states obey
\begin{dmath}\label{eq:universal-ground-state-identity}
	\Gamma(-\tau)\varphi_k(\tau)
	=
	\frac{1}{k}
	\ee^{-\pi\ii(k-1)\tau/k}
	\Gamma\!\left(-\frac{\tau}{k}\right)
\end{dmath}.
Introduce
\begin{dmath}\label{eq:Pearcey-Mellin-change}
	\tau_2 = {2t,
	\qquad
	\tau_4=-1-2t-2r
	}
\end{dmath}.
The factors depending on $\tau_2+\tau_4$ then become
\begin{dmath*}
	z^r
	\frac{\Gamma(1+r)\Gamma(-r)}
	{\Gamma(1+2r)}
\end{dmath*},
while the chosen branches give
\begin{dmath*}
	a_2^{2t}a_4^{\tau_4}
	\ee^{-\pi\ii t-3\pi\ii\tau_4/4}
	=
	2^{r+1/2}x^t
\end{dmath*}.
Including the vertical-variable Jacobian $4$ and the factors
$1/2$ and $1/4$ from
\cref{eq:universal-ground-state-identity}, we obtain
\begin{dmath*}
	\mathcal R(x,z)
	=
	\frac{1}{(2\pi\ii)^2}
	\int_{\mathcal C_t}\int_{\mathcal C_r}
	x^t z^r\Gamma(-t)\Gamma(-r)
	2^{r-1/2}
	\Gamma\!\left(\frac{2t+2r+1}{4}\right)
	\frac{\Gamma(1+r)}{\Gamma(1+2r)}
	\dd r\,\dd t
\end{dmath*}.
Finally,
\begin{dmath*}
	\frac{\Gamma(1+r)}{\Gamma(1+2r)}
	=
	\frac{\sqrt{\pi}\,2^{-2r}}{\Gamma(r+1/2)}
\end{dmath*},
so that
\begin{dmath}\label{eq:resolved-Pearcey-MB}
	\mathcal R(x,z)
	= {
	\frac{1}{(2\pi\ii)^2}
	\int_{\mathcal C_t}\int_{\mathcal C_r}
	x^t z^r\Gamma(-t)\Gamma(-r)\Phi(t,r)
	\dd r\,\dd t
	=
	G(x,z)
	}
\end{dmath}.
The transformed chamber is precisely
\begin{dmath*}
	\RE t< {0,
	\qquad
	\RE r<0,
	\qquad
	\RE(t+r)>-\frac12
	}
\end{dmath*}.
Thus the calculation recovers both the intrinsic ground state
\cref{eq:PE-ground-state} and the previously established
Pearcey Mellin--Barnes representation. Its vertical-decay
estimates also establish absolute convergence of the original
pairing \cref{eq:factorized-auf-resolvent} under the invertible
change of variables.

We have therefore proved, initially for $x,z>0$,
\begin{dmath}\label{eq:auf-interchange-pearcey}
	\left\langle
	\int_0^\infty
	2\ee^{-s\mathcal A(\boldsymbol u)}
	\cos(s\sqrt z)\dd s,
	\varphi_2(t_2)\varphi_4(t_4)
	\right\rangle
	=
	2\int_0^\infty
	\left\langle
	\ee^{-s\mathcal A(\boldsymbol u)},
	\varphi_2(t_2)\varphi_4(t_4)
	\right\rangle
	\cos(s\sqrt z)\dd s
	=
	G(x,z)
\end{dmath}.
The left-hand pairing is defined by the Mellin transform of
the integrated scalar kernel. Equality follows from two
independently defined analytic expressions, rather than from
an assumed interchange of improper integrals. It extends to
the common convergence sector and subsequently by analytic
continuation to the entire Euclidean field.

This is the additional content supplied by the AUF: it does not merely provide an additional evaluation rule for the rational symbol. It supplies the analytic realization from which such an extension is induced: the integrated scalar kernel is Mellin-transformable in a nonempty chamber, and its Mellin--Barnes pairing with the ground state reproduces the Pearcey field. The corresponding rational umbral expression may therefore be regarded as the formal image of this analytic construction.
The result is a pairing-level interchange theorem for
Euclidean Pearcey, not an operator-integral identity inside
the formal algebra. Other kernels require their own
Mellin-transformability, branch and convergence checks.

\section{Conclusions and outlook}
\label{sec:conclusions}

Higher-order multivariable Hermite--Kamp\'e de F\'eriet
polynomials provide the organizing principle of this work.
Their generating function expands the unfolding part of a
canonical catastrophe phase, while contour integration supplies
the moments of its leading exponential. On suitable decay
contours, the resulting Hermite expansion converges locally
uniformly in the unfolding parameters. It therefore provides
an analytic foundation for the subsequent umbral representations,
rather than merely a symbolic rearrangement.

The universal formal umbral operator $\mathfrak u$ expresses
these representations through explicit moment ground states.
A single umbra retains the polynomial dependence of the
unfolding, whereas several commuting umbrae assign one variable
to each control parameter and transfer their coupling to a
multivariate ground state. The fold, cusp and swallowtail
illustrate both descriptions. For the cusp, the Euclidean
formulation also connects the Pearcey field with quartic
scalar partition functions, including zero-dimensional
Higgs-type models.

The formal umbral framework is sufficient for these constructions when evaluation
precedes integration. Reversing the order asks for more:
integration of the linearized kernel produces a rational
expression whose evaluation is not determined by the original
formal rules. In the Pearcey example, its naive positive-power
evaluation vanishes by parity. Alternative prescriptions
require additional definitions and cannot be justified by
the formal exponential identity alone.

The specific added value of the AUF is to supply such an
extension through an analytic pairing. For Euclidean Pearcey,
the integrated scalar kernel has an explicit bivariate Mellin
transform. With specified branches and convergent contours,
its pairing with the monomial ground states reproduces the
independently established Pearcey representation. Integration
before pairing and pairing before integration thus give the
same field. This is a proved identity between analytic
pairings, not an assertion that an unevaluated formal operator
integral already exists. Nor does it imply that every formal
umbral expression automatically admits an AUF realization.

Several extensions merit investigation. Catastrophes involving
two integration variables, beginning with the umbilic integrals,
would test how effectively the Hermite construction handles
mixed polynomial terms. A related question is which features
of the umbral description reflect the underlying catastrophe
rather than the coordinates chosen to describe it. Natural
pairings between the associated Hermite polynomial families
may also clarify their algebraic and analytic relations.

For applications, a concrete next step is to use the Hermite
organization to generate correction terms in approximations
near caustics, and to test whether the AUF representations
facilitate their evaluation across different parameter regimes.
Any gain in accuracy or numerical stability must be established
by comparison with existing methods. The present result is
more specific: the Hermite construction yields exact formal
representations, while the AUF supplies, at least in the Euclidean
Pearcey case, a controlled extension that makes the reversed
order of integration and evaluation meaningful.

The main conclusion of this work is therefore not that formal umbral
calculus should be replaced by an analytic theory.  Rather, the
universal formal construction identifies precisely what the symbolic
operators mean and where their natural domain ends.  Within that domain
the formal method is exact and efficient.  When canonical integration
leads beyond entire functions of independent umbrae, the AUF provides a
controlled analytic extension at the level of Mellin pairings.  The
Pearcey integral supplies the first nontrivial example in which these
two levels can be compared explicitly.


\bibliographystyle{unsrtnat}
\bibliography{umbral_catastrophe_integrals.bib}

\end{document}